%% file: main.tex
\documentclass[10pt,twocolumn,letterpaper]{article}

\usepackage[pagenumbers]{cvpr} 

\makeatletter
\@namedef{ver@everyshi.sty}{}
\makeatother
\usepackage[mode=buildnew]{standalone}
\usepackage{tikz}
\usepackage{twoopt}
\usetikzlibrary{positioning, fit, backgrounds, arrows, calc, decorations.markings, arrows.meta, math, tikzmark, shapes, shapes.multipart}
\usepackage{makecell}
\usepackage{multirow}
\usepackage{xfp} 
\usepackage{transparent}
 \usepackage[accsupp]{axessibility}  

\usepackage{pgfplots}
\pgfplotsset{compat=1.11}
\usepackage{pgfplotstable}
\usepackage{longtable}
\pgfplotstableset{col sep=comma}
\pgfkeys{/pgf/number format/.cd, set thousands separator={}}
\pgfplotstableset{empty header/.style={every head row/.style={output empty row},}}
\newcommand{\insertTable}[3]{
\pgfplotstabletypeset[
    empty header,
    begin table=\begin{longtable},
    every first row/.append style={before row={%
    \caption{#2}%
    \label{#3}\\\toprule
    & LB $\uparrow$ & LB time [s] & UB $\downarrow$ & UB time [s] \\ \toprule    
    \endfirsthead
    \multicolumn{5}{c}%
    {{\bfseries Table \thetable\ Continued from previous page}} \\
    \toprule
    & LB $\uparrow$ & LB time [s] & UB $\downarrow$ & UB time [s] \\ \toprule 
    \endhead
    \midrule \multicolumn{5}{r}{\textit{Continued on next page}} \\ \bottomrule
    \endfoot
    \bottomrule
    \endlastfoot
    }},
    end table=\end{longtable},
    col sep=comma,
	columns/Instance/.style={string type, column type=l,string replace*={_}{-}},
	every row 0 column 0 Instance/.style={string replace*={MEAN}{\textsc{Mean}}},
	string replace*={_}{-},
	columns/LB/.style={fixed, precision=0},
	columns/UB/.style={fixed, precision=0},
 	columns/{LB time [s]}/.style={fixed, precision=1},
 	columns/{UB time [s]}/.style={fixed, precision=2},
	column type = r,
	every row no 0/.style={after row=\midrule},
]{#1}
}

\usepackage[shortlabels]{enumitem}

\usepackage{graphicx}
\usepackage{amsmath}
\usepackage{amssymb}
\usepackage{amsthm}
\usepackage{mathtools}
\usepackage{microtype}
\usepackage{booktabs} 
\usepackage{caption}
\usepackage{subcaption}
\usepackage{bbm} 
\usepackage{bbold}

\newcommand{\0}{\mathbb{0}}

\usepackage{twoopt}
\usepackage[noend,ruled,linesnumbered]{algorithm2e}
\DontPrintSemicolon 
\SetKwInput{KwInput}{Input}
\SetKwInput{KwOutput}{Output}

\usepackage[pagebackref,breaklinks,colorlinks]{hyperref}

\usepackage[capitalize]{cleveref}
\crefname{section}{Sec.}{Secs.}
\Crefname{section}{Section}{Sections}
\Crefname{table}{Table}{Tables}
\crefname{table}{Tab.}{Tabs.}

\DeclareMathOperator{\shp}{SP}
\DeclareMathOperator{\sign}{sign}
\newcommand{\pluseq}{\mathrel{+}=}
\newcommand{\minuseq}{\mathrel{-}=}

\providecommand{\abs}[1]{\lvert #1 \rvert}

\renewcommand{\P}{\mathcal{P}}

\providecommand{\X}{\mathcal{X}}
\providecommand{\I}{\mathcal{I}}
\providecommand{\J}{\mathcal{J}}

\providecommand{\R}{\mathbb{R}}

\theoremstyle{plain}

\newtheorem{prop}{Proposition}

\theoremstyle{definition}
\newtheorem{defn}{Definition}

\newtheorem{ex}{Example}

\theoremstyle{remark}
\newtheorem*{rem}{Remark}

\definecolor{customcolor1}{HTML}{800080}
\definecolor{customcolor2}{HTML}{990000}
\definecolor{customcolor3}{HTML}{009900}
\definecolor{customcolor4}{HTML}{808000}
\definecolor{customcolor5}{HTML}{D30E78}
\definecolor{customcolor6}{HTML}{036ffc}
\definecolor{customcolor7}{HTML}{eb9800}

\def\cvprPaperID{6021} 
\def\confName{CVPR}
\def\confYear{2022}

\begin{document}

\title{FastDOG: Fast Discrete Optimization on GPU}

\author{Ahmed Abbas \qquad Paul Swoboda \\
\small{Max Planck Institute for Informatics, Saarland Informatics Campus}}
\maketitle

\input{abstract}
\input{introduction}
\input{related-work}
\input{dual}
\input{primal}
\input{BDD}
\input{experiments}
\input{conclusion}

\input{acknowledgments}

{
\small
\bibliographystyle{ieee_fullname}
\bibliography{literature.bib}
}
\clearpage
\input{appendix}

\end{document}

%% file: abstract.tex
\begin{abstract}
We present a massively parallel Lagrange decomposition method for solving 0--1 integer linear programs occurring in structured prediction.
We propose a new iterative update scheme for solving the Lagrangean dual and a perturbation technique for decoding primal solutions.
For representing subproblems we follow~\cite{lange2021efficient} and use binary decision diagrams (BDDs).
Our primal and dual algorithms require little synchronization between subproblems and optimization over BDDs needs only elementary operations without complicated control flow.
This allows us to exploit the parallelism offered by GPUs for all components of our method.
We present experimental results on combinatorial problems from MAP inference for Markov Random Fields, quadratic assignment and cell tracking for developmental biology.
Our highly parallel GPU implementation improves upon the running times of the algorithms from~\cite{lange2021efficient} by up to an order of magnitude.
In particular, we come close to or outperform some state-of-the-art specialized heuristics while being problem agnostic. Our implementation is available at \url{https://github.com/LPMP/BDD}.
\end{abstract}

%% file: introduction.tex
\section{Introduction}
\label{sec:introduction}

Solving integer linear programs (ILP) efficiently on parallel computation devices is an open research question.
Done properly it would enable more practical usage of many ILP problems from structured prediction in computer vision and machine learning.
Currently, state-of-the-art generally applicable ILP solvers tend not to benefit much from parallelism~\cite{perumalla_gpu_mip_design_considerations}.
In particular, linear program (LP) solvers for computing relaxations benefit modestly (interior point) or not at all (simplex) from multi-core architectures.
In particular generally applicable solvers are not amenable for execution on GPUs.
To our knowledge there exists no practical and general GPU-based optimization routine and only a few solvers for narrow problem classes have been made GPU-compatible e.g.~\cite{abbas2021rama, shekhovtsov2016solving, tourani2018mplp++, xu2020fast_message_passing_mrf}. 
This, and the superlinear runtime complexity of general ILP solvers has hindered application of ILPs in large structured prediction problems, necessitating either restriction to at most medium problem sizes or difficult and time-consuming development of specialized solvers as observed for the special case of MAP-MRF~\cite{Kappes2015}.

\begin{figure}
\centering
\input{figures/solver_type_comparison}
\caption{Qualitative comparison of ILP solvers for structured prediction. Our solver (FastDOG) is faster than Gurobi~\cite{gurobi} and comparable to specialized CPU solvers, but outperformed by specialized GPU solvers. FastDOG is applicable to a diverse set of applications obviating the human effort for developing solvers for new problem classes.}
\label{fig:solver_tradeoff}
\end{figure}
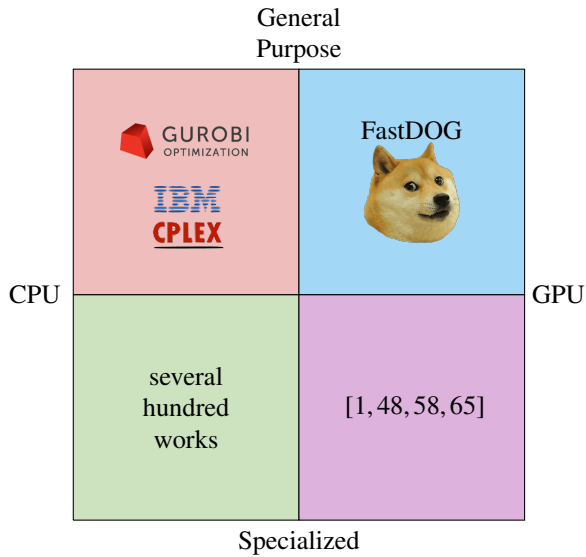

We argue that work on speeding up general purpose ILP solvers has had only limited success so far due to complicated control flow and computation interdependencies. 
We pursue an overall different approach and do not base our work on the typically used components of ILP solvers.
Our approach is designed from the outset to only use operations that offer sufficient parallelism for implementation on GPUs.

We argue that our approach sits on a sweet spot between general applicability and efficiency for problems in structured prediction as shown in Figure~\ref{fig:solver_tradeoff}.
Similar to general purpose ILP solvers~\cite{gurobi,cplex}, there is little or no effort to adapt these problems for solving them with our approach.
On the other hand we outperform general purpose ILP solvers in terms of execution speed for large problems from structured prediction and achieve runtimes comparable to hand-crafted specialized CPU solvers. We are only significantly outperformed by specialized GPU solvers.
However, development of fast specialized solvers especially on GPU is time-consuming and needs to be repeated for every new problem class.

Our work builds upon~\cite{lange2021efficient} in which the authors proposed a Lagrange decomposition into subproblems represented by binary decision diagrams (BDD).
The authors proposed sequential algorithms as well as parallel extensions for solving the Lagrange decomposition.
We improve upon their solver by proposing massively parallelizable GPU amenable routines for both dual optimization and primal rounding.
This results in significant runtime improvements as compared to their approach. 

%% file: figures/solver_type_comparison.tex
\begin{tikzpicture}[scale=3]

\definecolor{left_lower}{RGB}{205,227,191}
\definecolor{left_upper}{RGB}{239,189,187}
\definecolor{right_lower}{RGB}{221,179,221}
\definecolor{right_upper}{RGB}{162,215,245}

\draw[fill=left_lower] (0,0) -- (0,1) -- (1,1) -- (1,0) -- (0,0);
\draw[fill=left_upper] (0,1) -- (0,2) -- (1,2) -- (1,1) -- (0,1);
\draw[fill=right_lower] (1,0) -- (1,1) -- (2,1) -- (2,0) -- (1,0);
\draw[fill=right_upper] (1,1) -- (1,2) -- (2,2) -- (2,1) -- (1,1);


\node[align=center] (CPU) at (-0.17,1.0) {CPU}; 
\node[align=center] (GPU) at (2.15,1.0) {GPU}; 
\node (Specialized) at (1.0,-0.1) {Specialized};
\node[align=center] (General_Purpose) at (1.0,2.15) {General\\ Purpose};

\node[align=center] (Ours) at (1.5,1.5) 
{FastDOG\\ \includegraphics[width=50pt]{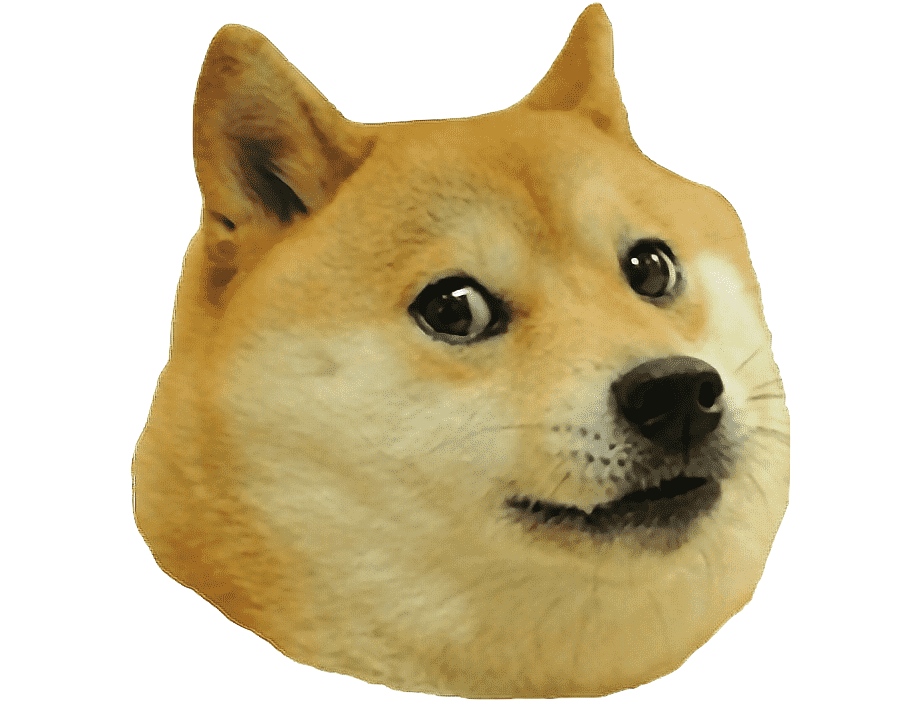}};
\node (spec_gpu) at (1.5,0.5) {
\hypersetup{citecolor=black} \cite{shekhovtsov2016solving,tourani2018mplp++,abbas2021rama,xu2020fast_message_passing_mrf}\hypersetup{citecolor=green}};
\node[align=center] (spec_gpu) at (0.5,0.5) {several\\ hundred\\works};
\node (Gurobi) at (0.5,1.68) {\includegraphics[width=50pt]{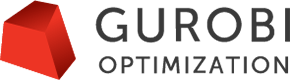}};
\node (ilog) at (0.5,1.32) {\includegraphics[width=50pt]{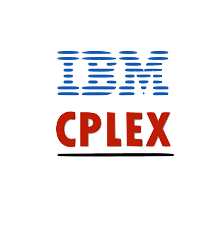}};

\end{tikzpicture}
 

%% file: related-work.tex
\section{Related Work}
\label{sec:related-work}

\paragraph{General Purpose ILP Solvers \& Parallelism}
The most efficient implementation of general purpose ILP solvers~\cite{gurobi,cplex} provided by commercial vendors typically benefit only moderately from parallelism. A recent survey is this direction is given in~\cite{perumalla_gpu_mip_design_considerations}. The main ways parallelism is utilized in ILP solvers are:
\begin{description}
\item[\textnormal{\textit{Multiple Independent Executions}}] State-of-the-art solvers~\cite{gurobi,cplex} offer the option of running multiple algorithms (dual/primal simplex, interior point, different parameters) solving the same problem in parallel until one finds a solution. While easy and worthwhile for problems for which best algorithms and parameters configurations are not known, such a simple approach can deliver parallelization speedups only to a limited degree.
\item[\textnormal{\textit{Parallel Branch-and-bound tree traversal}}] 
While appealing on first glance, it has been observed~\cite{ralphs2018parallel} that the order in which a branch-and-bound tree is traversed is crucial due to exploitation of improved lower and upper bounds and generated cuts. Consequently, it seems hard to obtain significant parallelization speedups and many recent improvements rely on a sequential execution. A separate line of work~\cite{sofranac2020accelerating_domain_propagation} exploited GPU parallelism for domain propagation allowing to decrease the size of the branch-and-bound tree.
\item[\textnormal{\textit{Parallel LP-Solver}}]
Interior point methods rely on computing a sequence of solutions of linear systems. This linear algebra can be parallelized for speeding up the optimization~\cite{gondzio2003parallel, smith12_gpu_interiorpoint}. However, for sparse problems sequential simplex solvers still outperform parallelized interior point methods. Also, a crossover step is needed to obtain a suitable basis for the simplex method for reoptimizing for primal rounding and in branch-and-bound searches, limiting the speedup obtainable by this sequential bottleneck. 
The simplex method is less straightforward to parallelize. The work~\cite{huangfu2018parallelizing} reports a parallel implementation, however current state-of-the-art commercial solvers outperform it with sequentially executed implementations.
\item[\textnormal{\textit{Machine Learning Methods}}]
Recently deep learning based methods have been proposed for choosing variables to branch on~\cite{gasse2019exact_co_gnn, nair2020solving} and for directly computing some easy to guess variables of a solution~\cite{nair2020solving} or improving a given one~\cite{sonnerat2021learning}.
While parallelism is not the goal of these works, the underlying deep networks are executed on GPUs and hence the overall computation heavy approach is fast and brings speedups.
Still, these parallel components do not replace the sequential parts of the solution process but work in conjunction with them, limiting the overall speedup attainable.
\end{description}
A shortcoming of the above methods in the application to very large structured prediction problems in machine learning and computer vision is that they still do not scale well enough to solve problems with  more than a millions variables in a few seconds.

\paragraph{Parallel Combinatorial Solvers}
For specialized combinatorial problem classes highly parallel algorithms for GPU have been developed.
For Maximum-A-Posteriori inference in Markov Random Fields~\cite{shekhovtsov2016solving, xu2020fast_message_passing_mrf} proposed a dual block coordinate ascent algorithm for sparse and~\cite{tourani2018mplp++} for dense graphs.
For multicut a primal-dual algorithm has been proposed in~\cite{abbas2021rama}.
Max-flow GPU implementations have been investigated in~\cite{vineet2008cuda,wu2012chapter}.
While some parts of the above specialized algorithms can potentially be generalized, other key components cannot, limiting their applicability to new problem classes and requiring time-consuming design of algorithms whenever attempting to solve a different problem class.
 
\paragraph{Specialized CPU solvers}
There is a large literature of specialized CPU solvers for specific problem classes in structured prediction.
For an overview of pursued algorithmic techniques for the special case of MRFs we refer to the overview article~\cite{Kappes2015}.
Most related to our approach are the so called dual block coordinate ascent (a.k.a.\ message passing) algorithms which optimize a Lagrange decomposition.
Solvers have been developed for MRFs~\cite{kolmogorov2006convergent,kolmogorov2014new,globerson2008fixing,werner2007linear,savchynskyy12efficient,jancsary2011convergent,meltzer2012convergent,wang2013subproblem,johnson2007lagrangian,tourani2018mplp++,tourani2020taxonomy}, graph matching~\cite{zhang2016pairwise,swoboda2017study,swoboda2019convex}, multicut~\cite{lange2018partial,swoboda2017message,abbas2021rama}, multiple object tracking~\cite{hornakova2021making} and cell tracking~\cite{haller2020primal}.
Most of the above algorithms require a sequential computation of update steps.

\paragraph{Optimization with Binary Decision Diagrams} 
Our work builds upon~\cite{lange2021efficient}.
The authors proposed a Lagrange decomposition of ILPs that can be optimized via a sequential dual block coordinate ascent method or a decomposition based approach that can utilize multiple CPU cores.

The works~\cite{BergmanCire2016,bergman2018discrete,lozano2018consistent} similarly consider decompositions into multiple BDDs and solve the resulting problem with general purpose ILP solvers.
The work~\cite{bergman2015lagrangian} investigates optimization of Lagrange decompositions with multi-valued decision diagrams with subgradient methods.
An extension for job sequencing was proposed in~\cite{improved_job_sequence_bounds_hooker_2019} and in~\cite{castro2020mdd} for routing problems.
Hybrid solvers using mixed integer programming solvers were investigated in~\cite{tjandraatmadja2020incorporating,gonzalez2020bdd,gonzalez2020integrated}.
The works~\cite{andersen2007constraint,bergman2016decision,bergman2016discrete} consider stable set and max-cut and propose optimizing (i)~a relaxation to get lower bounds~\cite{andersen2007constraint} or (ii)~a restriction to generate approximate solutions~\cite{bergman2016decision,bergman2016discrete}.

In contrast to previous BDD-based optimization methods we propose a highly parallelizable and problem agnostic approach that is amenable to GPU computation.

%% file: dual.tex
\section{Method}
\label{sec:problem}

\begin{table}
    \centering
    \begin{tabular}{ll}
         \toprule
         $x_i$ & Optimization variable $i \in [n]$\\
         $\X_j$ & Feasible set of constraint $j \in [m]$ \\
         $\I_j$ & Set of variables in constraint $j \in [m]$\\
         $\J_i$ & Set of constraints containing variable $i \in [n]$\\
         $m_{ij}^{\beta}$ & \makecell[l]{Min-marginal for variable $i$ taking value $\beta$ in \\ subproblem $j \in [m]$} \\
         $\lambda_i^j$ & Lagrange multiplier for variable $i$ in subproblem $j$ \\
         \bottomrule
         \end{tabular}
    \caption{Notation of symbols used in our problem decomposition.}
    \label{table:notation}
\end{table}

\begin{figure*}
\centering
\includestandalone[width=\textwidth]{figures/overview_figure}
\caption{Example decomposition of a binary program into two subproblems, one for each constraint. Each subproblem is represented by a weighted BDD where solid arcs model the cost $\lambda$ of assigning a $1$ to the variable and dashed arcs have $0$ cost which model assigning a $0$. All $r-\top$ paths in BDDs encode feasible variable assignments of corresponding subproblems (and $r-\bot$ infeasible). Optimal assignments w.r.t current (non-optimal) $\lambda$ are highlighted in green i.e.~$x_1 = 1, x_2 = x_3 = 0$ for $\X_1$ and $x_2=x_3=x_4=0$ for $\X_2$. Our dual update scheme processes multiple variables in parallel which are indicated in same color (e.g. $x_1, x_2$ in $\X_1, \X_2$ resp.).}
\label{fig:overview}
\end{figure*}

We first introduce the optimization problem and its Lagrange decomposition.
Next we elaborate our parallel update scheme for optimizing the Lagrangean dual followed by our parallel primal rounding algorithm.
For the problem decomposition and dualization we follow~\cite{lange2021efficient}. 
Our notation is summarized for reference in Table~\ref{table:notation}.

\begin{defn}[Binary Program]
Consider a linear objective $c \in \R^n$ and
$m$ variable subsets $\I_j \subset [n]$ of constraints with feasible set $\X_j \subset \{0,1\}^{\I_j}$ for $j \in [m]$.
The corresponding binary program is defined as
\begin{equation}
\label{eq:binary-program}
    \min_{x \in \{0,1\}^n} c^\top x \quad \text{s.t.} \quad x_{\I_j} \in \X_j \quad \forall j \in [m]\,,
\tag{BP} 
\end{equation}
where $x_{\I_j}$ is the restriction to variables in $\I_j$.
\end{defn}

\begin{ex}[ILP]
    \label{ex:ILP}
 Consider the 0--1 integer linear program
\begin{align}
\min \quad & c^\top x \quad
\text{s.t.} \quad Ax \leq b, \; x \in \{0,1\}^n.
\tag{ILP}
\label{eq:ILP}
\end{align}
The system of linear constraints $Ax \leq b$ may be split into $m$ blocks, each block representing a single (or multiple) rows of the system. For instance, let $a_j^\top x \leq b_j$ denote the $j$-th row of $Ax \leq b$, then the problem can be written in the form \eqref{eq:binary-program} by setting $\I_j = \{ i \in [n] : a_{ji} \neq 0\}$ and $\X_j = \{ x \in \{0,1\}^{\I_j} : \sum_{i \in \I_j} a_{ji} x_{i} \leq b_j \}$.
\end{ex}

\subsection{Lagrangean Dual}

While~\eqref{eq:binary-program} is NP-hard to solve, optimization over a single constraint is typically easier for example by using Binary Decision Diagrams. To make use of this possibility we dualize the original problem using Lagrange decomposition similarly to~\cite{lange2021efficient}.
This allows us to solve the Lagrangean dual of the full problem~\eqref{eq:binary-program} by solving only the subproblems.

\begin{defn}[Lagrangean dual problem]
Define the set of subproblems that constrain variable $x_i$ as
$\J_i = \{j \in [m] \mid i \in \I_j\}$.
Let the energy for subproblem $j \in [m]$ w.r.t.\ Lagrangean dual variables $\lambda^j \in \R^{\I_j}$ be
\begin{equation}
E^j(\lambda^j) = \min_{x \in \X_j} x^\top \lambda^{j} \,.
\end{equation}
Then the Lagrangean dual problem is defined as
\begin{align}
\max_\lambda \quad \sum_{j \in [m]} E^j(\lambda^j) \quad
\text{s.t.} \quad \sum_{j \in \J_i} \lambda^j_i = c_i \quad \forall i \in [n]. \label{eq:dual-problem} \tag{D}
\end{align}
\end{defn}

If optima of the individual subproblems $E^j(\lambda^j$) agree with each other then the consensus vector obtained from stitching together individual subproblem solutions solves the original problem~\eqref{eq:binary-program}.
In general,~\eqref{eq:dual-problem} is a lower bound on~\eqref{eq:binary-program}.
Formal derivation of~\eqref{eq:dual-problem} is given in~\cite{lange2021efficient}.

\subsection{Min-Marginals}
To optimize the dual problem~\eqref{eq:dual-problem} and also to obtain a primal solution we use min-marginals~\cite{lange2021efficient} defined as

\begin{defn}[Min-marginals]
For $i \in [n]$, $j \in \J_i$ and $\beta \in \{0,1\}$ let
\begin{align}
m^\beta_{ij} = \min_{x \in \X_j} x^\top \lambda^{j} \quad \text{s.t.} \quad x_i = \beta \label{eq:min-marginals}
\tag{MM}
\end{align}
denote the \emph{min-marginal} w.r.t.\ primal variable $i$, subproblem $j$ and $\beta$.
\end{defn}
\begin{defn}[Min-marginal differences]
For notational convenience let us also define 
\begin{align}
M_{ij} = m_{ij}^1 - m_{ij}^0, 
\label{eq:min-marginals-diff}
\tag{MD}
\end{align}
which denote \emph{min-marginal difference} computed through~\eqref{eq:min-marginals}.
\end{defn}
If $M_{ij} > 0$ then assigning a value of $0$ to variable $i$ has a lower cost than assigning a $1$ in the subproblem $j$ and viceversa. Thus, the quantity $\abs{M_{ij}}$ indicates by how much $E^j(\lambda^j)$ increases if $x_i$ is fixed to $1$ (if $M_{ij} > 0$), respectively $0$ (if $M_{ij} < 0$). 

Min-marginals have been used in various ways to design dual block coordinate ascent algorithms~\cite{lange2021efficient,arora2013higher,zhang2016pairwise,swoboda2017study,swoboda2019convex,swoboda2017message,haller2020primal,kolmogorov2006convergent,kolmogorov2014new,globerson2008fixing,werner2007linear,savchynskyy12efficient,jancsary2011convergent,meltzer2012convergent,wang2013subproblem,johnson2007lagrangian,tourani2018mplp++,tourani2020taxonomy,werner2019relative,abbas2021rama,hornakova2021making}.

\begin{algorithm}
\KwInput{Lagrange variables $\lambda_i^j \in \R \, \forall i \in [n], j \in \J_i$,
Constraint sets $\X_j \subset \{0,1\}^{\I_j} \, \forall j \in [m]$, 
Damping factor $\omega \in (0,1]$}
Initialize deferred min-marginal diff.  $\overline{M} = \0$\;
\label{alg:initialization-line}
\While{(stopping criterion not met)}
{
\For{$j \in \J$ in parallel}
{
\label{alg:variable-order-loop-start}
\For{$i \in \I_j$ in ascending order}
{
\label{alg:min-marginal-computation-start}
Compute min-marginal diff. $M_{ij}$ \eqref{eq:min-marginals-diff}\;
\label{alg:min-marginal-computation-end}%
Update dual variables $\lambda_i^j$ via \eqref{eq:dual_update_parallel}
\label{alg:lambda-update}
} 
}
Update deferred min-marginal diff. $\overline{M} \leftarrow M$\;
\label{alg:deferred-min-marginals-update}
Repeat lines~\ref{alg:variable-order-loop-start}-\ref{alg:lambda-update} in descending order of $\I_j$\;
\label{alg:optimization-end}
}
\For{$j \in \J$, $i \in \I_j$}{
Add deferred min-marginal differences:
$\lambda^{j}_i \pluseq \omega \overline{M}_{ij}$\;
\label{alg:lambda-last-update}
}
\caption{Parallel Deferred Min-Marg. Averaging}
\label{alg:parallel-mma}
\end{algorithm}

\subsection{Parallel Deferred Min-Marginal Averaging}
To exploit GPU parallelism in solving the dual problem~\eqref{eq:dual-problem} we would like to update multiple dual variables in parallel. However, conventional dual update schemes are not friendly for parallelization. For example the dual update scheme of~\cite{lange2021efficient} for variable $i$ in subproblem $j$ is
\begin{equation}
    \lambda_i^{j} \leftarrow \lambda_i^{j} - M_{ij} + \underbrace{\frac{1}{\abs{\J_i}} \sum_{k \in \J_i} M_{ik}\,,}_{\text{min-marginal averaging}}
    \label{eq:dual_update_serial}
\end{equation}
where $M_{ij}$ is defined in~\eqref{eq:min-marginals-diff}.
This update scheme~\eqref{eq:dual_update_serial} requires communication between all subproblems $\J_i$ containing variable $i$ for the min-marginal averaging step and thus requires synchronization.
To overcome this limitation we propose a novel dual optimization procedure which performs this averaging step on min-marginal differences $\overline{M}$ from the previous iteration as follows
\begin{equation}
    \lambda_i^{j} \leftarrow \lambda_i^{j} - \omega M_{ij} + \frac{\omega}{\abs{\J_i}} \sum_{k \in \J_i} \overline{M}_{ik}.
    \label{eq:dual_update_parallel}
\end{equation}
Since $\overline{M}$ was computed in the previous iteration, the above dual updates can be performed in parallel for all subproblems without requiring synchronization. Following~\cite{werner2019relative} we use a damping factor $\omega \in (0,1)$ ($0.5$ in our experiments) to obtain better final solutions. 

Our proposed scheme is given in Algorithm~\ref{alg:parallel-mma}. We iterate in parallel over each subproblem $j$. For each subproblem, variables are visited in order and min-marginals are computed and stored for updates in the next iteration (lines~\ref{alg:min-marginal-computation-start}-\ref{alg:min-marginal-computation-end}).
The current min-marginal difference is subtracted and the one from previous iteration is added (line~\ref{alg:lambda-update}) by distributing it equally among subproblems $\J_i$. At termination (line~\ref{alg:lambda-last-update}) we perform a min-marginal averaging step to account for the deferred update from last iteration. For stopping criteria we use relative change in dual objective between two subsequent iterations. We initialize the input Lagrange variables by $\lambda^{j}_i = c_i / \abs{\J_i}$, $\forall i \in [n]$, $j \in \J_i$.
\begin{prop}
\label{prop:convergence}
In each dual iteration the Lagrange multipliers along with the deferred min-marginals can be used to satisfy dual feasibility and the dual lower bound~\eqref{eq:dual-problem} is non-decreasing.
\end{prop}
Similar to other dual block coordinate ascent schemes Algorithm~\ref{alg:parallel-mma} can get stuck in suboptimal points, see~\cite{werner2019relative,werner2007linear}.
As seen in our experiments these are usually not far away from the optimum, however.

In Section~\ref{sec:bdd} we will explain how we can incrementally compute min-marginals reusing previous computations if we represent subproblems as Binary Decision Diagrams.
This saves us from computing min-marginals from scratch leading to greater efficiency. 

%% file: figures/overview_figure.tex
\begin{tikzpicture}[align=center, x=0cm, y=0cm, >=stealth, line width=0.15mm,
myrect/.style={draw, rectangle, rounded corners, inner sep=1mm},
every node/.style={inner sep = 1mm, minimum size=3mm, align=center, font=\scriptsize},
mycircle/.style={draw, circle, inner sep = 0.5mm},
node distance = 0.7cm]
\tikzmath{\columnsep = 130;} 

\node[myrect, font = \scriptsize] (lp) {$\min_{x \in \{0,1\}} -5x_1 + x_2 + 4x_3 + 3x_4$ \\
$\textcolor{black}{\X_1}: x_1 + x_2 + x_3 \leq 2$, \quad \qquad $\textcolor{black}{\X_2}: x_2 + x_3 - x_4 = 0$};

\node[mycircle, fill = customcolor1!40] at (-\columnsep, -2cm) (r1) {\tiny $r$};
\node[draw=none, label=left:{\textcolor{gray}{$x_1$}},fit=(r1)] (l1) {};

\node[mycircle, fill = customcolor6!30, above right = of r1] (n1_a) {};
\node[mycircle, fill = customcolor6!30, below right = of r1] (n1_b) {};
\node[myrect, gray, densely dotted, label=above:{\textcolor{gray}{$x_2$}},fit=(n1_a) (n1_b)] (l2) {};

\node[mycircle, fill = customcolor4!30, right = of n1_a] (n2_a) {};
\node[mycircle, fill = customcolor4!30, right = of n1_b] (n2_b) {};
\node[myrect, gray, densely dotted, label=above:{\textcolor{gray}{$x_3$}},fit=(n2_a) (n2_b)] (l3) {};


\node[mycircle,  right = of n2_a] (n3_a) {\tiny $\top$};
\node[mycircle,  right = of n2_b] (n3_b) {\tiny $\bot$};

\draw[->, dashed] (r1) to (n1_a);
\draw[->, customcolor3!50] (r1) to node [draw=none, below left, pos=0.5, font=\tiny] {$-5$} (n1_b);

\draw[->] (n1_a) to [out = 22, in = 157, looseness = 1] node [draw=none, above, midway, font=\tiny] {$0.5$} (n2_a);
\draw[->, dashed] (n1_a) to [out = -22, in = 202, looseness = 1] (n2_a);

\draw[->, dashed, customcolor3!50] (n1_b) to (n2_a);
\draw[->] (n1_b) to node [draw=none, below, midway, font=\tiny] {$0.5$} (n2_b);

\draw[->] (n2_a) to [out = 22, in = 157, looseness = 1] node [draw=none, above, midway, font=\tiny] {$2$} (n3_a);
\draw[->, dashed, customcolor3!50] (n2_a) to [out = -22, in = 202, looseness = 1] (n3_a);

\draw[->, dashed] (n2_b) to (n3_a);
\draw[->] (n2_b) to node [draw=none, below, midway, font=\tiny] {$2$} (n3_b);

\node[draw=none, below right = 0.5cm and -2.7cm of n2_b, font = \scriptsize] (bdd2_label) {
$\lambda^1 = (-5, 0.5, 2)$, $(x_1, x_2, x_3) \in \textcolor{black}{\X_1}$ };




\node[mycircle, fill = customcolor1!40] at (\columnsep - 90, -2cm) (r2) {\tiny $r$};
\node[draw=none, label=left:\textcolor{gray}{$x_2$},fit=(r2)] (l1_2) {};

\node[mycircle, fill = customcolor6!30, above right = of r2] (m1_a) {};
\node[mycircle, fill = customcolor6!30, below right = of r2] (m1_b) {};
\node[myrect, gray, densely dotted, label=above:\textcolor{gray}{$x_3$},fit=(m1_a) (m1_b)] (l2_2) {};

\node[mycircle, fill = customcolor4!30, right = of m1_a] (m2_a) {};
\node[mycircle, fill = customcolor4!30, right = of m1_b] (m2_b) {};
\node[myrect, gray, densely dotted, label=above:\textcolor{gray}{$x_4$},fit=(m2_a) (m2_b)] (l3_2) {};

\node[mycircle, right = of m2_a] (m3_a) {\tiny $\top$};
\node[mycircle, right = of m2_b] (m3_b) {\tiny $\bot$};

\draw[->] (r2) to node [draw=none, inner sep = 1mm, below, pos=0.3, font=\tiny] {$0.5$\quad} (m1_b);
\draw[->, dashed, customcolor3!50] (r2) to (m1_a);

\draw[->] (m1_b) to [out = -45, in = 225, looseness = 1] node [draw=none, below, pos=0.125, font=\tiny] {$2$} (m3_b);

\draw[->] (m1_a) to node [draw=none, above, pos=0.5, font=\tiny] {$2$} (m2_b);

\draw[->, dashed, customcolor3!50] (m1_a) to (m2_a);
\draw[->, dashed] (m1_b) to (m2_b);

\draw[->, dashed, customcolor3!50] (m2_a) to (m3_a);
\draw[->, dashed] (m2_b) to (m3_b);

\draw[->] (m2_a) to node [draw=none, above, pos=0.3, font=\tiny] {\quad $3$} (m3_b);

\draw[->] (m2_b) to node [draw=none, below, pos=0.3, font=\tiny] { \quad $3$} (m3_a);
\node[draw=none, below right = 0.5cm and -1.8cm of m1_b, font=\scriptsize] (bdd2_label) {$\lambda^2 = (0.5, 2, 3)$, $(x_2, x_3, x_4) \in \textcolor{black}{\X_2}$};

\end{tikzpicture}

%% file: primal.tex
\section{Primal Rounding}
\label{sec:primal}


In order to obtain a primal solution to~\eqref{eq:binary-program} from an approximative dual solution to~\eqref{eq:dual-problem} we propose a GPU friendly primal rounding scheme based on cost perturbation. We iteratively change costs in a way that variable assignments across subproblems agree with each other. If all variables agree by favoring a single assignment, we can reconstruct a primal solution (not necessarily the optimal). Instead of only using variable assignments of all subproblems we use min-marginal differences~\eqref{eq:min-marginals-diff} as they additionally indicate how strongly a variable favours a particular assignment. 

Algorithm~\ref{alg:primal-rounding} details our method. We iterate over all variables in parallel and check min-marginal differences.
If for a variable $i$ all min-marginal differences indicate that the optimal solution is $0$ (resp.\ 1) Lagrange variables $\lambda$ are increased (resp.\ decreased) leaving even more certain min-marginals differences for these variables. This step imitates variable fixation as done in branch-and-bound, however we only perform soft fixation implicitly through cost perturbation. 
In case min-marginal differences are equal we randomly perturb corresponding dual costs. Lastly, if min-marginals differences indicate conflicting solutions we compute total min-marginal difference and decide accordingly. In the last two cases we add more perturbation to force towards non-conflicts. For faster convergence we increase the perturbation magnitude after each iteration. 

Note that the modified $\lambda$ variables via Alg.~\ref{alg:primal-rounding} need not be feasible for the dual problem~\eqref{eq:dual-problem}. Although, our primal rounding algorithm is not guaranteed to terminate, in our experiments a solution was always found in less than $100$ iterations. 

\begin{algorithm}
\KwInput{
Lagrange variables $\lambda_i^j \in \R \, \forall i \in [n], j \in \J_i$,
Constraint sets $\X_j \subset \{0,1\}^{\I_j} \, \forall j \in [m]$,
Initial perturbation strength $\delta \in \R_+$, 
perturbation growth rate $\alpha$
}
\KwOutput{Feasible labeling $x \in \{0,1\}^n$}
Compute min-marginal differences $M_{ij}\, \forall i, j$ ~\eqref{eq:min-marginals-diff}\;
\While{$\exists i \in [n]$ and $j\neq k \in \J_i$ s.t.\ $\sign(M_{ij}) \neq \sign(M_{ik})$} {
\For{$i=1,\ldots,n$ in parallel } {
Sample $r$ uniformly from $[-\delta,\delta]$\;
\If{$M_{ij} > 0$ $\forall j \in \J_i$} {
$\lambda_i^{j} \pluseq \delta$ \quad $\forall j \in \J_i$\; 
}
\ElseIf{$M_{ij} < 0$ $\forall j \in \J_i$} {
$\lambda_i^{j} \minuseq \delta$ \quad $\forall j \in \J_i$\; 
}
\ElseIf{$ M_{ij} = 0 $ $\forall j \in \J_i$} {
$\lambda_i^{j} \pluseq r \cdot \delta$ \quad $\forall j \in \J_i$\; 
}
\Else
{
Compute total min-marginal difference: $M_i = \sum_{j \in \J_i} M_{ij}$\;
    $\lambda_i^{j} \pluseq \sign(M_i) \cdot \lvert r \vert \cdot \delta$ \quad $\forall j \in \J_i$\; 
}
}
Increase perturbation: $\delta \leftarrow \delta \cdot \alpha$\;
Reoptimize perturbed $\lambda$ via Algorithm~\ref{alg:parallel-mma}\;
Recompute $M_{ij}\, \forall i, j$ w.r.t optimized $\lambda$ \;
}
\caption{Perturbation Primal Rounding}
\label{alg:primal-rounding}
\end{algorithm}

\begin{rem}
The primal rounding scheme in~\cite{lange2021efficient} and typical primal ILP heuristics~\cite{berthold2006primal} are sequential and build upon sequential operations such as variable propagation.
Our primal rounding lends itself to parallelism since we perturb costs on all variables simultaneously and reoptimize via Algorithm~\ref{alg:parallel-mma}.
\end{rem}

%% file: BDD.tex
\section{Binary Decision Diagrams}
\label{sec:bdd}

We use Binary Decision Diagrams (BDDs) to represent the feasibility sets $\X_j$, $j \in [m]$ and compute their min-marginals~\eqref{eq:min-marginals}.
BDDs are in essence directed acyclic graphs whose paths between two special nodes (root and terminal) encode all feasible solutions. Specifically, we use reduced ordered Binary Decision Diagrams~\cite{bryant1986graph} as in~\cite{lange2021efficient}.

\begin{defn}[BDD]
Let an ordered variable set $\I = \{w_1,\ldots,w_k\} \subset [n]$ corresponding to a constraint be given.
A corresponding BDD is a directed acyclic graph $D=(V,A)$ with
\begin{description}
\item[\normalfont\textit{Special nodes}:]
 root node $r$, terminals $\bot$ and $\top$.
\item[\normalfont\textit{Outgoing Arcs}:]
each node $v \in V \backslash \{\top, \bot\}$ has exactly two successors $s^0(v), s^1(v)$ with outgoing arcs $v s^0(v) \in A$ (the zero arc) and $v s^1(v) \in A$ (the one arc).
\item[\normalfont\textit{Partition}:]
the node set $V$ is partitioned by $\{\P_1, \ldots, \P_k\}$, $\dot{\cup}_i \P_i = V \backslash \{ \top, \bot\}$.
Each partition holds all the nodes corresponding to a single variable \eg $\P_i$ corresponds to variable $w_i$.
It holds that $\P_1= \{r\}$ \ie it only contains the root node.
\item[\normalfont\textit{Partition Ordering}:]
when $v \in \P_i$ then $s^0(v), s^1(v) \in \P_{i+1} \cup \{\bot\}$ for $i < k$ and $s^0(v),s^1(v) \in \{\bot,\top\}$ for $v \in \P_k$. 
\end{description}
\label{def:bdd}
\end{defn}

\begin{defn}[Constraint Set Correspondence] 
Each BDD defines a constraint set $\X$ via the relation
\begin{equation}
    x \in \X \Leftrightarrow 
    \begin{array}{c}
    \exists (v_1,\ldots,v_{k},v_{k+1}) \in \text{Paths}(V,A) \text{ s.t.\ } \\ v_1 = r, v_{k+1} = \top, \\
    v_{i+1} = s^{x_i}(v_i)\, \forall i \in [k]
    \end{array}
\end{equation}
Thus each path between root $r$ and terminal $\top$ in the BDD corresponds to some feasible variable assignment $x \in \X$.
\label{def:bdd_correspondence}
\end{defn}
Figures~\ref{fig:overview} and ~\ref{fig:weighted_BDD} illustrate BDD encoding of feasible sets of linear inequalities.

\begin{rem}
In the literature~\cite{bryant1986graph,knuth2011art} BDDs have additional requirements, mainly that there are no isomorphic subgraphs.
This allows for some additional canonicity properties like uniqueness and minimality.
While all the BDDs in our algorithms satisfy the additional canonicity properties, only what is required in Definition~\ref{def:bdd} is needed for our purposes, so we keep this simpler setting.
\end{rem}

\subsection{Efficient Min-Marginal Computation}
In order to compute min-marginals for subproblems we need to consider weighted BDDs. For notational convenience we will drop dependence on subproblem $j$ in upcoming text \eg, we will use $\lambda_{i}$ instead of $\lambda_{i}^j$. 

\begin{defn}[Weighted BDD]
A weighted BDD is a BDD with arc costs.
Let a function $f(x)$ be defined as
\begin{equation}
    f(x) = 
    \begin{cases}
        x^\top \lambda & x \in \X\\
        \infty & \text{otherwise}
    \end{cases}.
\end{equation}
The weighted BDD represents $f$ if it satisfies Def.~\ref{def:bdd_correspondence} for the given $\X$ and the arc costs for an $i \in [k], v \in \P_i,  vw \in A$ are set as
    $\begin{cases}
    0 & w \in s^0(v) \\
    \lambda_i & w \in s^1(v) \\
    \end{cases}$.
\end{defn}

\begin{figure}
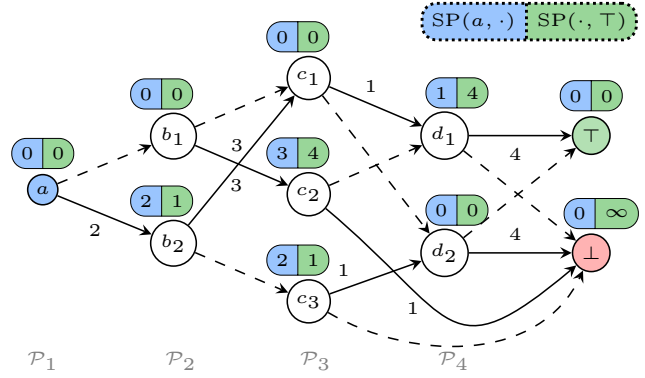

\centering
\includestandalone[width=\columnwidth]{figures/bdd_ops}
\caption{Weighted BDD of a subproblem containing variables: $\I = \{a, b, c, d\}$ with costs ($\lambda$): $2, 3, 1, 4$ resp.~and constraint $a - b - c + d = 0$. Shortest path costs from the root node $a$ and the terminal node $\top$ are shown for each node. Here $\P_1 = \{a\}, \P_2 = \{b_1, b_2\}, \P_3 = \{c_1, c_2, c_3\}, \P_4 = \{d_1, d_2\}$, $s^0(c_2) = d_1$ and $s^1(c_2) = \bot$. Dashed arcs have cost $0$ as they model assigning a $0$ value to the corresponding variable. }
\label{fig:weighted_BDD}
\end{figure}

Min-marginals for variable $i \in \I$ of a subproblem can be computed by its weighted BDD by calculating shortest path distances from $r$ to all nodes in $\P_i$ and shortest path distances from all nodes in $\P_{i+1}$ to $\top$.
We use $\shp(v,w)$ to denote the shortest path distance between nodes $v$ and $w$ of a weighted BDD. An example shortest path calculation is shown in Figure~\ref{fig:weighted_BDD}.
The min-marginals as defined in~\eqref{eq:min-marginals} can be computed as
\begin{equation}
\label{eq:min-marginal-via-shortest-path}
   m^{\beta}_i = \min_{\substack{vs^{\beta}(v) \in A\\ v \in \P_i}} \left[\shp(r,v) + \beta\cdot\lambda_i + \shp(s^{\beta}(v),\top)\right]
\end{equation}

For efficient min-marginal computation in Algorithm~\ref{alg:parallel-mma} we reuse shortest path distances used in~\eqref{eq:min-marginal-via-shortest-path}.
Specifically, for computing min-marginals in  lines~\ref{alg:min-marginal-computation-start}-\ref{alg:min-marginal-computation-end} of Alg.~\ref{alg:parallel-mma} we use Alg.~\ref{alg:forward_pass_mm} for ascending variable order and Alg.~\ref{alg:backward_pass_mm} for descending variable order in line~\ref{alg:optimization-end}. 
\begin{algorithm}
    \For{$v \in \P_i$}{
        $\shp(r,v) = \min\left\{
        \begin{array}{c}
            \min\limits_{u: s^0(u) = v}\ \shp(r,u),\\
            \min\limits_{u: s^1(u) = v}\ \shp(r,u) + \lambda_i 
        \end{array}
        \right\}$
        }
        Compute $m_i^{\beta}$ via~\eqref{eq:min-marginal-via-shortest-path}\;
    \caption{Forward Pass Min-Marginal Computation}
    \label{alg:forward_pass_mm}
\end{algorithm}

\begin{algorithm}
    \For{$v \in \P_{i+1}$}{
        $\shp(v,\top) = \min\left\{
        \begin{array}{c}
            \shp(s^0(v),\top), \\ \shp(s^1(v),\top) + \lambda_{i+1}
        \end{array}
        \right\}$        
        }
    Compute $m_i^{\beta}$ via~\eqref{eq:min-marginal-via-shortest-path}\;
    \caption{Backward Pass Min-Marginal Computation}
    \label{alg:backward_pass_mm}
\end{algorithm}

\paragraph{Efficient GPU implementation}
In addition to solving all subproblems in parallel, we also exploit parallelism within each subproblem during shortest path updates. Specifically in Alg.~\ref{alg:forward_pass_mm}, we parallelize over all $v \in \P_i$ and perform the $\min$ operation atomically. Similarly in Alg.~\ref{alg:backward_pass_mm} we parallelize over all $v \in \P_{i+1}$ but without requiring atomic update. 

\begin{figure}
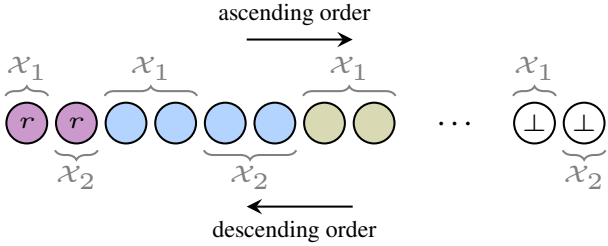

\centering
\includestandalone[width=\columnwidth]{figures/bdd_memory_arrangement}
\caption{Arrangement of BDD nodes in GPU memory for the ILP in Figure~\ref{fig:overview}. For ascending order in Alg.~\ref{alg:parallel-mma} we proceed from root to terminal nodes and vice versa for descending.}
\label{fig:memory_arrangement}
\end{figure}

To enable fast GPU memory access via memory coalescing we arrange BDD nodes in the following fashion. First, all nodes within a BDD which belong to the same partition $\P$ (thus corresponding to same variable) are laid out consecutively. Secondly, across different BDDs, nodes are ordered w.r.t increasing hop distance from their corresponding root nodes. Such arrangement for the ILP in Figure~\ref{fig:overview} is shown in Figure~\ref{fig:memory_arrangement}.

%% file: figures/bdd_ops.tex
\begin{tikzpicture}[align=center, x=0cm, y=0cm, >=stealth, line width=0.15mm,
myrect/.style={draw, rectangle, rounded corners, inner sep=1mm},
every node/.style={inner sep = 1mm, minimum size=3mm, align=center},
mycircle/.style={draw, circle, inner sep = 0.5mm},
state/.style={
rectangle split,
rectangle split horizontal,
rectangle split parts=2,
rectangle split part fill={customcolor6!40, customcolor3!40},
rounded corners,
draw=black, ultra thin,
inner sep=0.9mm,
text centered},
node distance = 0.78cm]
\tikzmath{\columnsep = 130;} 
\tikzstyle{every node}=[font=\fontsize{4.5}{1}\selectfont]

\node[mycircle, fill = customcolor6!40] at (0, 0) (r1) {$a$};

\node[mycircle, above right = 0.25cm and 1cm of r1] (n1_a) {$b_1$};
\node[mycircle, below right = 0.25cm and 1cm of r1] (n1_b) {$b_2$};

\node[mycircle, above right = 0.25cm and 1cm of n1_a] (n2_a) {$c_1$};
\node[mycircle, below right = 0.25cm and 1cm of n1_a] (n2_b) {$c_2$};
\node[mycircle, below right = 0.25cm and 1cm of n1_b] (n2_c) {$c_3$};

\node[mycircle, above right = 0.25cm and 1cm of n2_b] (n3_a) {$d_1$};
\node[mycircle, below right = 0.25cm and 1cm of n2_b] (n3_b) {$d_2$};

\node[mycircle, fill = customcolor3!30, right = 1cm of n3_a] (n4_a) {\tiny $\top$};
\node[mycircle, fill = red!30, right = 1cm of n3_b] (n4_b) {\tiny $\bot$};

\draw[->, dashed] (r1) to (n1_a);
\draw[->] (r1) to node [draw=none, below, pos=0.4, font=\tiny] {$2$} (n1_b);

\draw[->, dashed] (n1_a) to (n2_a);
\draw[->] (n1_a) to node [draw=none, above, pos=0.45, font=\tiny] {$3$} (n2_b);

\draw[->, dashed] (n1_b) to (n2_c);
\draw[->] (n1_b) to node [draw=none, below, pos=0.45, font=\tiny] {$3$} (n2_a);

\draw[->, dashed] (n2_a) to (n3_b);
\draw[->] (n2_a) to node [draw=none, above, pos=0.45, font=\tiny] {$1$} (n3_a);
\draw[->, dashed] (n2_b) to (n3_a);
\draw[->] (n2_c) to node [draw=none, above, pos=0.15, font=\tiny] {$1$} (n3_b);
\draw[->] (n2_b) to [out = -40, in = 220, looseness = 1.7] node [draw=none, below, pos=0.3, font=\tiny] {$1$} (n4_b);
\draw[->, dashed] (n2_c) to [out = -30, in = 245, looseness = 1] (n4_b);

\draw[->] (n3_a) to node [draw=none, below, pos=0.45, font=\tiny] {$4$} (n4_a);
\draw[->, dashed] (n3_a) to (n4_b);
\draw[->, dashed] (n3_b) to (n4_a);
\draw[->] (n3_b) to node [draw=none, fill= white, above, pos=0.45, font=\tiny] {$4$} (n4_b);

\node[state, above = 0.05cm of r1] (r1_c) {$0$ \nodepart{two} $0$};
\node[state, above left = 0.1cm and -0.35cm of n1_a] (n1_a_c) {$0$ \nodepart{two} $0$};
\node[state, above left = 0.1cm and -0.35cm of n1_b] (n1_b_c) {$2$ \nodepart{two} $1$};

\node[state, above left = 0.1cm and -0.35cm of n2_a] (n2_a_c) {$0$ \nodepart{two}$0$};
\node[state, above left = 0.1cm and -0.35cm of n2_b] (n2_b_c) {$3$ \nodepart{two}$4$};
\node[state, above left = 0.1cm and -0.35cm of n2_c] (n2_c_c) {$2$ \nodepart{two}$1$};

\node[state, above right = 0.1cm and -0.35cm of n3_a] (n3_a_c) {$1$ \nodepart{two} $4$};
\node[state, above right = 0.1cm and -0.35cm of n3_b] (n3_b_c) {$0$ \nodepart{two} $0$};

\node[state, above = 0.05cm of n4_a] (n4_a_a) {$0$ \nodepart{two} $0$};
\node[state, above right = 0.1cm and -0.4cm of n4_b] (n4_b_a) {$0$ \nodepart{two} $\infty$};

\node[state, thick, densely dotted, above right = 0.75cm and -0.4cm of n3_a] (ann) {$\shp(a, \cdot)$ \nodepart{two} $\shp(\cdot, \top)$};

\node[draw=none, gray, below = 1.3cm of r1] (r1t) {$\P_1$};
\node[draw=none, gray, right = of r1t] (n1t) {$\P_2$};
\node[draw=none, gray, right = 0.75cm of n1t] (n2t) {$\P_3$};
\node[draw=none, gray, right = of n2t] (n3t) {$\P_4$};

\end{tikzpicture}

%% file: figures/bdd_memory_arrangement.tex
\begin{tikzpicture}[align=center, x=0cm, y=0cm, >=stealth, line width=0.15mm,
myrect/.style={draw, rectangle, rounded corners, inner sep=1mm},
every node/.style={inner sep = 1mm, minimum size=3mm, align=center, font=\scriptsize},
mycircle/.style={draw, circle, inner sep = 0.1mm},
node distance = 0.05cm]

\node[mycircle, fill = customcolor1!40] (r1) {\tiny $r$};
\node[mycircle, fill = customcolor1!40, right = of r1] (r2) {\tiny $r$};

\node[mycircle, fill = customcolor6!30, right = of r2] (n1) {};
\node[mycircle, fill = customcolor6!30, right = of n1] (n2) {};
\node[mycircle, fill = customcolor6!30, right = of n2] (n3) {};
\node[mycircle, fill = customcolor6!30, right = of n3] (n4) {};

\node[mycircle, fill = customcolor4!30, right = of n4] (n5) {};
\node[mycircle, fill = customcolor4!30, right = of n5] (n6) {};

\node[draw=none, right = 0.5cm of n6] (n10) {};
\node[mycircle, right = of n10] (n11) {\tiny $\bot$};
\node[mycircle, right = of n11] (n12) {\tiny $\bot$};

\node at ($(n6)!.75!(n10)$) {\ldots};

\node[draw=none, above = 0.3cm of n3] (n4a) {};
\node[draw=none, above = 0.3cm of n6] (n6a) {};
\draw[->] (n4a) to node [draw=none, above, pos=0.45, font=\tiny] {ascending order} (n6a);

\node[draw=none, below = 0.3cm of n3] (n4b) {};
\node[draw=none, below = 0.3cm of n6] (n6b) {};
\draw[<-] (n4b) to node [draw=none, below, pos=0.45, font=\tiny] {descending order} (n6b);

\draw [decorate, gray, decoration={brace, amplitude=2pt, raise=2mm}] (r1.west) -- (r1.east) node [black, midway, yshift=4mm, gray]{\tiny$\X_1$};

\draw [decorate, gray, decoration={brace, amplitude=2pt, raise=2mm}] (n1.west) -- (n2.east) node [black, midway, yshift=4mm, gray]{\tiny$\X_1$};

\draw [decorate, gray, decoration={brace, amplitude=2pt, raise=2mm}] (n5.west) -- (n6.east) node [black, midway, yshift=4mm, gray]{\tiny$\X_1$};

\draw [decorate, gray, decoration={brace, mirror, amplitude=2pt, raise=2mm}] (r2.west) -- (r2.east) node [black, midway, yshift=-4mm, gray]{\tiny$\X_2$};

\draw [decorate, gray, decoration={brace, mirror, amplitude=2pt, raise=2mm}] (n3.west) -- (n4.east) node [black, midway, yshift=-4mm, gray]{\tiny$\X_2$};

\draw [decorate, gray, decoration={brace, amplitude=2pt, raise=2mm}] (n11.west) -- (n11.east) node [black, midway, yshift=4mm, gray]{\tiny$\X_1$};

\draw [decorate, gray, decoration={brace, mirror, amplitude=2pt, raise=2mm}] (n12.west) -- (n12.east) node [black, midway, yshift=-4mm, gray]{\tiny$\X_2$};

\end{tikzpicture}

%% file: experiments.tex
\section{Experiments}
\label{sec:experiments}

\begin{table*}
    \centering
    \resizebox{\textwidth}{!}{
    \begin{tabular}{l rr rrr rrrr rr}
    \toprule
    \multirow{2}{*}{} & \multicolumn{2}{c}{\textit{Cell tracking}} & \multicolumn{3}{c}{\textit{Graph matching}} &  \multicolumn{4}{c}{\textit{MRF}} & \multicolumn{2}{c}{\textit{QAPLib}} \\
    \cmidrule(lr){2-3} \cmidrule(lr){4-6} \cmidrule(lr){7-10} \cmidrule(lr){11-12}
     & \textit{Small} & \textit{Large} & \textit{Hotel} & \textit{House} &  \textit{Worms} & \textit{C-seg} & \textit{C-seg-n4} & \textit{C-seg-n8} & \textit{Obj-seg} &  \textit{Small} & \textit{Large} \\     
    \midrule
    \# instances & 10 & 5 & 105 & 105 & 30 & 3 & 9 & 9 & 5 & 105 & 29 \\
    $n_{max}$ & $1.2M$ & $10M$ & $0.3M$ & $0.3M$ & $1.5M$ & $3.3M$ & $1.2M$ & $1.4M$ & $681k$ & $3M$ & $49M$ \\
    $m_{max}$ & $0.2M$ & $2.3M$ & $52k$ & $52k$ & $0.2M$ & $13.6M$ & $4.2M$ & $8.3M$ & $2.2M$ & $245k$ & $2M$ \\
    \midrule
    \multicolumn{12}{c}{Dual objective (lower bound) $\uparrow$} \\
    \midrule
    \texttt{Gurobi}~\cite{gurobi} & 
    $-\textbf{4.382e6}$ & $-\textbf{1.545e8}$ & 
    $-\textbf{4.293e3}$ & $-\textbf{3.778e3}$ & $-4.849e4$ & 
    $\textbf{3.085e8}$ & $1.9757e4$ & $1.9729e4$ & $3.1311e4$ &
    $2.913e6$ & $4.512e4$     \\ 
    \texttt{BDD-CPU}~\cite{lange2021efficient} &
    $-4.387e6$ & $-1.549e8$ &
    $-\textbf{4.293e3}$ & $-\textbf{3.778e3}$ & $-4.878e4$ &
    ${3.085e8}$ & $1.9643e4$ & $1.9631e4$ & $3.1248e4$ &
    $3.675e6$ & $8.172e6$ \\
    \texttt{Specialized} &
    $-4.385e6$ & $-1.551e8$ & 
    $-\textbf{4.293e3}$ &  $-\textbf{3.778e3}$ &  $-\textbf{4.847e4}$    &  
    $\textbf{\textbf{3.085e8}}$ & $\textbf{2.0012e4}$ &  $\textbf{1.9991e4}$ & $\textbf{3.1317e4}$  & 
    - & -    \\ 
    \texttt{FastDOG} &
    $-4.387e6$ & $-1.549e8$ & 
    $-\textbf{4.293e3}$    &   $-\textbf{3.778e3}$   &  $-4.893e4$ &  
    ${3.085e8}$ & ${2.0011e4}$ &   ${1.9990e4}$ &  $\textbf{3.1317e4}$   & 
    $\textbf{3.747e6}$ &  \textbf{8.924e6}   \\ 
    \midrule
    \multicolumn{12}{c}{Primal objective (upper bound) $\downarrow$} \\
    \midrule
    \texttt{Gurobi}~\cite{gurobi} & 
    $-\textbf{4.382e6}$ & $-1.524e8$ & $-\textbf{4.293e3}$ & $-\textbf{3.778e3}$ & $-4.842e4$ & $\textbf{3.085e8}$ & $2.8464e4$ & $2.7829e4$ & $1.4981e5$ & $5.186e7$ & $1.431e8$ \\
    \texttt{BDD-CPU}~\cite{lange2021efficient} &
    $-4.337e6$ & $-1.515e8$ & $-\textbf{4.293e3}$ & $-\textbf{3.778e3}$ & $-4.783e4$ & $3.086e8$ & $2.1781e4$ & $2.2338e4$ & $3.1525e4$ & $5.239e7$ & $1.452e8$ \\
    \texttt{Specialized} &
     $-4.361e6$ & $-1.531e8$ & 
    $-\textbf{4.293e3}$ & $-\textbf{3.778e3}$ &   $-\textbf{4.845e4}$   &  
    $\textbf{3.085e8}$ & $\textbf{2.0012e4}$ &  $\textbf{1.9991e4}$ & $\textbf{3.1317e4}$  & 
      - & -    \\ 
          
    \texttt{FastDOG} &
    $-4.376e6$  & $-\textbf{1.541e8}$ & 
    $-\textbf{4.293e3}$ &   $-\textbf{3.778e3}$   &  $-4.831e4$    &  
    ${3.085e8}$   &  $2.0016e4$  &   ${1.9995e4}$ & ${3.1322e4}$ & 
    $\textbf{4.330e7}$  &  $\textbf{1.376e8}$   \\ 
    \midrule
    \multicolumn{12}{c}{Runtimes [s] $\downarrow$} \\
    \midrule
    \texttt{Gurobi}~\cite{gurobi}  &
    $\textbf{1}$ & $1584$ & $4$ & $7$ & $1048$ & $132$ & $980$ & $1337$ & $1506$ & $3948$ & $6742$ \\
    \texttt{BDD-CPU}~\cite{lange2021efficient} &
    $14$ & $216$ & $6$ & $12$ & $528$ & $70$ & $107$ & $218$ & $232$ & $357$ & $\textbf{5952}$ \\
    \texttt{Specialized} &
    $1.5$ & $\textbf{90}$ & 
     ${3}$   &    ${3}$  & $214$ &  
    $155$ &  $\textbf{9}$ &    ${30}$  & $\textbf{3}$ & 
         - & -     \\ 
    \texttt{FastDOG}     &
    $\fpeval{3 + 10}$ & $\fpeval{39 + 71} $  & 
    $\textbf{0.2}$   & $\textbf{0.4}$ &   $\textbf{54}$   &  
    $\textbf{14}$ &  $\textbf{9}$   &  $\textbf{13}$ &  $39$ & 
    $137$     &   \fpeval{2713+4215}  \\ 
    \bottomrule
    \end{tabular}}
    \caption{Results comparison on all datasets where the values are averaged within a dataset. For each dataset, the results on corresponding specialized solvers are computed using~\cite{haller2020primal, swoboda2017study, kappes_trws_pct}. Numbers in bold highlight the best performance. $n_{max}, m_{max}$: Maximum number of variables, constraints in the category. }
    \label{tab:combined-results}
\end{table*}

\begin{figure*}
\captionsetup[subfigure]{justification=centering}
\centering
     \begin{subfigure}[b]{0.33\textwidth}
         \includegraphics[width=\textwidth,height=5cm]{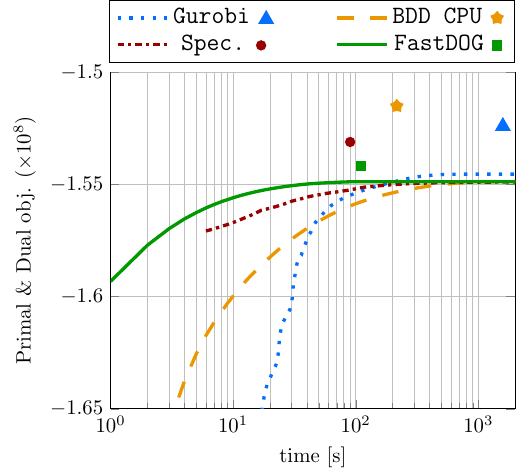}
         \caption{\textit{Cell tracking: Large}}
     \end{subfigure}
     \begin{subfigure}[b]{0.33\textwidth}
         \includegraphics[width=\textwidth,height=5cm]{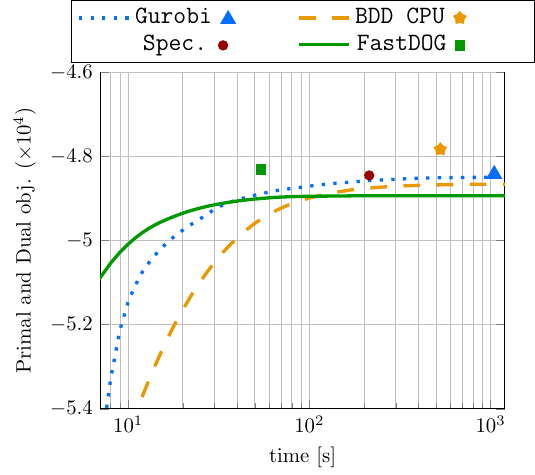}
         \caption{\textit{Graph Matching: Worms}}
     \end{subfigure}
     \begin{subfigure}[b]{0.33\textwidth}
         \includegraphics[width=\textwidth,height=5cm]{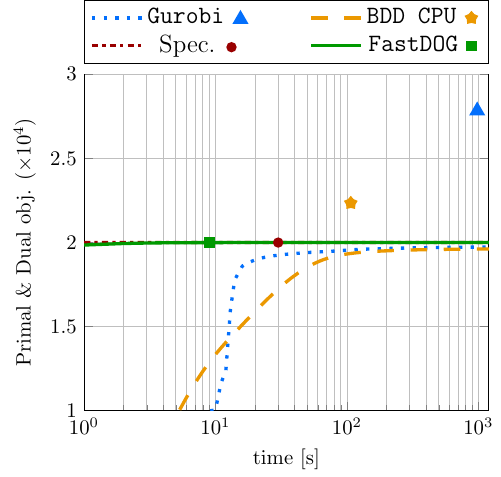}
         \caption{\textit{MRF: Color-seg-n8}}
     \end{subfigure}
    \caption{Convergence plots averaged over all instances of a dataset. Lower curves depict increasing lower bounds while markers denote objectives of rounded primal solutions. The x-axis is plotted logarithmically.}
    \label{fig:convergence-plots}
\end{figure*}

We show effectiveness of our solver against a state-of-the-art ILP solver~\cite{gurobi}, the general purpose BDD-based solver~\cite{lange2021efficient} and specialized CPU solvers for specific problem classes.
We have chosen some of the largest structured prediction ILPs we are aware of in the literature that are publicly available. Our results are computed on a single NVIDIA Volta V100 (16GB) GPU unless stated otherwise. For CPU solvers we use AMD EPYC 7702 CPU. 

\paragraph{Datasets}
Our benchmark problems obtained from~\cite{swoboda2022structured_pred_archive} can be categorized as follows.

\begin{description}[itemsep=2pt,parsep=2pt]
\item[\textnormal{\textit{Cell tracking}}:] Instances from~\cite{haller2020primal} which we partition into small and large instances as also done in~\cite{lange2021efficient}.

\item[\textnormal{\textit{Graph matching (GM)}}:] Quadratic assignment problems (often called graph matching in the literature) for correspondence in computer vision~\cite{torresani2008feature} (\texttt{hotel}, \texttt{house}) and developmental biology~\cite{kainmueller2014active} (\texttt{worms}).

\item[\textnormal{\textit{Markov Random Field (MRF)}}:] Several datasets from the OpenGM~\cite{Kappes2015} benchmark, containing both small and large instances with varying topologies and number of labels. We have chosen the datasets \textit{color-seg}, \textit{color-seg-n4}, \textit{color-seg-n8} and \textit{object-seg}.
\item[\textnormal{\textit{QAPLib}}:] The widely used benchmark dataset for quadratic assignment problems used in the combinatorial optimization community~\cite{QAPLIB}.
We partition QAPLib instances into small (up to $50$ vertices) and large (up to $128$ vertices) instances. For large instances we use NVIDIA RTX 8000 (48GB) GPU.
\end{description}
\paragraph{Algorithms}
We compare results on the following algorithms.
\begin{description}[itemsep=2pt,parsep=2pt]
	\item[\normalfont\texttt{Gurobi}:] The commercial ILP solver~\cite{gurobi} as reported in~\cite{lange2021efficient}.
	The barrier method is used for QAPLib and dual simplex for all other datasets.
	\item[\normalfont\texttt{BDD-CPU}:] BDD-based min-marginal averaging approach of~\cite{lange2021efficient}.
	The algorithm runs on CPU with 16 threads for parallelization.
	Primal solutions are rounded using their BDD-based depth-first search scheme.
	\item[\normalfont\texttt{Specialized solvers}:] State-of-the-art problem specific solver for each dataset.
	For cell-tracking we use the solver from~\cite{haller2020primal}, the \texttt{AMP} solver for graph matching proposed in~\cite{swoboda2017study} and \texttt{TRWS} for MRF~\cite{kolmogorov2006convergent}.
	\item[\normalfont\texttt{FastDOG}:] Our approach where for the GPU implementation we use the CUDA~\cite{cuda} and Thrust~\cite{Thrust} programming frameworks.
	For rounding primal solutions with Algorithm~\ref{alg:primal-rounding} we set $\delta = 1.0$ and $\alpha = 1.2$.
	For constructing BDDs out of linear (in)equalities we use the same approach as for \texttt{BDD-CPU}.
\end{description}
For \textit{MRF}, parallel algorithms such as~\cite{tourani2018mplp++} exist however \texttt{TRWS} is faster on the sparse problems we consider.
While we are aware of even faster purely primal heuristics~\cite{komodakis2007approximate,boykov2001fast} for \textit{MRF} and \eg~\cite{hutschenreiter2021fusion} for \textit{graph matching} they do not optimize a convex relaxation and hence do not provide lower bounds.
Hence, we have chosen \texttt{TRWS}~\cite{kolmogorov2006convergent} for \textit{MRF} and \texttt{AMP}~\cite{swoboda2017study} for \textit{graph matching} which, similar to \texttt{FastDOG}, optimize an equivalent resp.\ similar Lagrange decomposition and hence can be directly compared. 

\paragraph{Results}
In Table~\ref{tab:combined-results} we show aggregated results over all instances of each specific benchmark dataset.
Runtimes are taken w.r.t.\ computation of both primal and dual bounds.
A more detailed table with results for each instance is given in the Appendix.

In Figure~\ref{fig:convergence-plots} we show averaged convergence plots for various solvers.
In general we offer a very good anytime performance producing at most times and in general during the beginning better lower bounds than our baselines.

\paragraph{Discussion}
In general, we are always faster (up to a factor of $10$) than \texttt{BDD-CPU}~\cite{lange2021efficient} and except on \textit{worms} we achieve similar or better lower bounds.
In comparison to the respective hand-crafted \texttt{Specialized} CPU solvers we also achieve comparable runtimes with comparable lower and upper bounds.
While \texttt{Gurobi} achieves, if given unlimited time, better lower bounds and primal solutions, our \texttt{FastDOG} solver outperforms it on the larger instances when we abort \texttt{Gurobi} early after hitting a time limit.
We argue that we outperform \texttt{Gurobi} on larger instances due to its superlinear iteration complexity.

When comparing the number of dual iterations to \texttt{BDD-CPU} we need roughly $3$-times as many to reach the same lower bound. 
Nonetheless, as we can perform more iterations per second this still leads to an overall faster algorithm.

Since we are solving a relaxation the lower bounds and quality of primal solutions are dependent on the tightness of this relaxation. For all datasets except \textit{QAPLib} our (and also baselines') lower and upper bounds are fairly close, reflecting the nature of commonly occurring structured prediction problems.

%% file: conclusion.tex
\section{Conclusion}
\label{sec:conclusion}
We have proposed a massively parallelizable generic algorithm that can solve a wide variety of ILPs on GPU.
Our results indicate that the performance of specialized efficient CPU solvers can be matched or even surpassed by a completely generic GPU solver.
Our implementation is a first prototype and we conjecture that more speedups can be gained by elaborate implementation techniques, \eg compression of the BDD representation, better memory layout for better memory coalescing, multi-GPU support etc.
We argue that future improvements in optimization algorithms for structured prediction can be made by developing GPU friendly problem specific solvers and with improvements in our or other generic GPU solvers that can benefit many problem classes simultaneously. 
Another future avenue is optimization of ILPs from other domains, \eg on the MIPLib benchmark~\cite{miplib_2017}.
These problems include constraints that are harder to represent as BDDs and additional encoding techniques are needed~\cite{abio2012new,bdd_ordering_fujita}.

%% file: acknowledgments.tex
\section*{Acknowledgments}
We would like to thank all anonymous reviewers for their valuable feedback and especially Reviewer 1 for very detailed reading and identifying areas of improvement. 
We also thank Jan-Hendrik Lange for insightful discussions.

%% file: appendix.tex
\onecolumn
\section*{Appendix}
\label{sec:appendix}

\section{Proof of Proposition~\ref{prop:convergence}}
\begin{proof}
\noindent\textbf{\\ Feasibility of iterates}
We prove 
\begin{equation}
\sum_{j \in \J_i} \lambda^{j}_i + \omega M_{ij} = c_i
\end{equation}
just before line~\ref{alg:deferred-min-marginals-update} in Algorithm~\ref{alg:parallel-mma}.
We do an inductive proof over the number of iterates w.r.t iterations $t$.
\begin{description}
\item[$t = 0$:] Follows from $\overline{M} = 0$ and the uniform distribution of costs in line~\ref{alg:initialization-line} of Algorithm~\ref{alg:parallel-mma}.
\item[$t > 0$:] Let $\lambda(t-1)$, $M(t-1)$, $\overline{M}(t-1)$ be last iterations' Lagrange multipliers, min-marginals differences and (deferred) min-marginal differences. Also let $\lambda(t)$, $M(t)$ and $\overline{M}(t)$ be the same from current iteration just before line~\ref{alg:deferred-min-marginals-update}. Note that $\overline{M}(t) = M(t-1)$.
It holds that
\begin{subequations}
\begin{align}
\sum_{j \in \J_i} \left[\lambda^{j}_i(t) + \omega M_{ij}(t)\right]
=&\sum_{j \in \J_i} \left[\lambda^j_i(t-1) - \omega M_{ij}(t) + \sum_{k \in \J_i} \left(\frac{\omega}{\abs{\J_i}} \overline{M}_{ik}(t) \right)  + \omega M_{ij}(t) \right] \\
=& \sum_{j \in \J_i} \left[\lambda^j_i(t-1) + \omega M_{ij}(t-1)\right] \\
=& c_i\,.
\end{align}
\end{subequations}

\end{description}

\noindent\textbf{Non-decreasing Lower Bound}
In order to prove that iterates have non-decreasing lower bound we will consider an equivalent lifted representation in which proving the non-decreasing lower bound will be easier.

\noindent\textbf{\\ Lifted Representation}
Introduce $\lambda^{j,\beta}_i$ for $\beta \in \{0,1\}$ and the subproblems
\begin{equation}
E(\lambda^{j,1}, \lambda^{j,0}) = \min_{x \in \X_j} x^\top \lambda^{j,1} + (1-x)^\top \lambda^{j,0}
\end{equation}
Then~\eqref{eq:dual-problem} is equivalent to
\begin{equation}
\label{eq:lifted-dual-problem}
\max_{\lambda^{1},\lambda^0} \sum_{j \in \J} E(\lambda^{j,1}, \lambda^{j,0}) \text{ s.t. } \sum_{j \in \J_i} \lambda^{j,\beta}_i = \beta \cdot c_i
\end{equation}
We have the transformation from original to lifted $\lambda$
\begin{equation}
\lambda \mapsto (\lambda^1 \leftarrow \lambda, \lambda^0 \leftarrow \0)
\end{equation}
and from lifted to original $\lambda$ (except a constant term)
\begin{equation}
(\lambda^1, \lambda^0) \mapsto \lambda^1 - \lambda^0\,.
\end{equation}
It can be easily shown that the lower bounds are invariant under the above mappings and feasible $\lambda$ for~\eqref{eq:dual-problem} are mapped to feasible ones for~\eqref{eq:lifted-dual-problem} and vice versa. 

The update rule line~\ref{alg:lambda-update} in Algorithm~\ref{alg:parallel-mma} for the lifted representation can be written as
\begin{equation}
\label{eq:lifted-representation-lambda-update}
\lambda_i^{j,\beta} \leftarrow \lambda_i^{j,\beta} - \omega \cdot \min(m^{\beta}_{ij} - m^{1-\beta}_{ij},0) + \omega \cdot \min(\overline{m}^{\beta}_{ij} - \overline{m}^{1-\beta}_{ij},0)
\end{equation}
It can be easily shown that~\eqref{eq:lifted-representation-lambda-update} and line~\ref{alg:lambda-update} in Algorithm~\ref{alg:parallel-mma} are corresponding to each other under the transformation from lifted to original $\lambda$.

\noindent\textbf{\\ Continuation of Non-decreasing Lower Bound}
Define
\begin{equation}
\lambda_i^{\prime j,\beta} = \lambda_i^j - \omega \cdot \min(m^{\beta}_{ij} - m^{1-\beta}_{ij},0)\,.
\end{equation}
Then $E(\lambda^{\prime j,1}, \lambda^{\prime j,0}) = E(\lambda^{j,1}, \lambda^{j,0})$ 
are equal due to the definition of min-marginals.
Define next
\begin{equation}
\lambda_i^{\prime \prime j,\beta} = \lambda_i^{\prime j} + \omega \cdot \min(\overline{m}^{\beta}_{ij} - \overline{m}^{1-\beta}_{ij},0)\,.
\end{equation}
Then $E(\lambda^{\prime \prime j,1}, \lambda^{\prime \prime j,0}) \geq E(\lambda^{\prime j,1}, \lambda^{\prime j,0})$ since $\lambda'' \geq \lambda'$. 
This proves the claim. 
Another side-benefit of the lifted update scheme~\eqref{eq:lifted-representation-lambda-update} is that evaluating $E(\lambda^{j,1}, \lambda^{j,0})$ during the course of optimization in Algorithm~\ref{alg:parallel-mma} always gives a lower bound to the true dual objective calculated after accounting for deferred min-marginal differences.
\end{proof}




\section{Detailed results}
\subsection{\textit{Cell tracking}}
\insertTable{figures/tables/ct-small.csv}{Detailed results of \texttt{FastDOG} on \textit{Cell tracking - small} dataset}{tab:ct-small-results}
\insertTable{figures/tables/ct-large.csv}{Detailed results of \texttt{FastDOG} on \textit{Cell tracking - large} dataset}{tab:ct-large-results}

\subsection{\textit{Graph matching}}
\insertTable{figures/tables/gm-hotel.csv}{Detailed results of \texttt{FastDOG} on \textit{Graph matching - hotel} dataset}{tab:gm-hotel-results}
\insertTable{figures/tables/gm-house.csv}{Detailed results of \texttt{FastDOG} on \textit{Graph matching - house} dataset}{tab:gm-house-results}
\insertTable{figures/tables/gm-worms.csv}{Detailed results of \texttt{FastDOG} on \textit{Graph matching - worms} dataset}{tab:gm-worms-results}

\subsection{\textit{MRF}}
\insertTable{figures/tables/color-seg.csv}{Detailed results of \texttt{FastDOG} on \textit{MRF - color-seg} dataset}{tab:mrf-color-seg-results}
\insertTable{figures/tables/color-seg-n4.csv}{Detailed results of \texttt{FastDOG} on \textit{MRF - color-seg-n4} dataset}{tab:mrf-color-seg-n4-results}
\insertTable{figures/tables/color-seg-n8.csv}{Detailed results of \texttt{FastDOG} on \textit{MRF - color-seg-n8} dataset}{tab:mrf-color-seg-n8-results}
\insertTable{figures/tables/object-seg.csv}{Detailed results of \texttt{FastDOG} on \textit{MRF - object-seg} dataset}{tab:mrf-object-seg-results}

\subsection{\textit{QAPLib}}
\insertTable{figures/tables/qaplib-small.csv}{Detailed results of \texttt{FastDOG} on \textit{QAPLib - small} dataset}{tab:qaplib-small}
\insertTable{figures/tables/qaplib-large.csv}{Detailed results of \texttt{FastDOG} on \textit{QAPLib - large} dataset}{tab:qaplib-large}
\clearpage
\twocolumn

%% file: main.bbl
\begin{thebibliography}{10}\itemsep=-1pt

\bibitem{abbas2021rama}
Ahmed Abbas and Paul Swoboda.
\newblock {RAMA: A Rapid Multicut Algorithm on GPU}.
\newblock {\em arXiv preprint arXiv:2109.01838}, 2021.

\bibitem{abio2012new}
Ignasi Ab{\'\i}o, Robert Nieuwenhuis, Albert Oliveras, Enric
  Rodr{\'\i}guez-Carbonell, and Valentin Mayer-Eichberger.
\newblock A new look at bdds for pseudo-boolean constraints.
\newblock {\em Journal of Artificial Intelligence Research}, 45:443--480, 2012.

\bibitem{andersen2007constraint}
Henrik~Reif Andersen, Tarik Hadzic, John~N Hooker, and Peter Tiedemann.
\newblock A constraint store based on multivalued decision diagrams.
\newblock In {\em International Conference on Principles and Practice of
  Constraint Programming}, pages 118--132. Springer, 2007.

\bibitem{arora2013higher}
Chetan Arora and Amir Globerson.
\newblock Higher order matching for consistent multiple target tracking.
\newblock In {\em Proceedings of the IEEE International Conference on Computer
  Vision}, pages 177--184, 2013.

\bibitem{BergmanCire2016}
David Bergman and Andre~A. Cire.
\newblock Decomposition based on decision diagrams.
\newblock In Claude-Guy Quimper, editor, {\em Integration of AI and OR
  Techniques in Constraint Programming}, pages 45--54, Cham, 2016. Springer
  International Publishing.

\bibitem{bergman2018discrete}
David Bergman and Andre~A Cire.
\newblock Discrete nonlinear optimization by state-space decompositions.
\newblock {\em Management Science}, 64(10):4700--4720, 2018.

\bibitem{bergman2015lagrangian}
David Bergman, Andre~A Cire, and Willem-Jan van Hoeve.
\newblock Lagrangian bounds from decision diagrams.
\newblock {\em Constraints}, 20(3):346--361, 2015.

\bibitem{bergman2016decision}
David Bergman, Andre~A Cire, Willem-Jan Van~Hoeve, and John Hooker.
\newblock {\em Decision diagrams for optimization}, volume~1.
\newblock Springer, 2016.

\bibitem{bergman2016discrete}
David Bergman, Andre~A Cire, Willem-Jan van Hoeve, and John~N Hooker.
\newblock Discrete optimization with decision diagrams.
\newblock {\em INFORMS Journal on Computing}, 28(1):47--66, 2016.

\bibitem{berthold2006primal}
Timo Berthold.
\newblock Primal heuristics for mixed integer programs.
\newblock 2006.

\bibitem{boykov2001fast}
Yuri Boykov, Olga Veksler, and Ramin Zabih.
\newblock Fast approximate energy minimization via graph cuts.
\newblock {\em IEEE Transactions on pattern analysis and machine intelligence},
  23(11):1222--1239, 2001.

\bibitem{bryant1986graph}
Randal~E Bryant.
\newblock Graph-based algorithms for boolean function manipulation.
\newblock {\em Computers, IEEE Transactions on}, 100(8):677--691, 1986.

\bibitem{QAPLIB}
Rainer~E Burkard, Stefan~E Karisch, and Franz Rendl.
\newblock {QAPLIB}--a quadratic assignment problem library.
\newblock {\em Journal of Global optimization}, 10(4):391--403, 1997.

\bibitem{castro2020mdd}
Margarita~P Castro, Andre~A Cire, and J~Christopher Beck.
\newblock An mdd-based lagrangian approach to the multicommodity
  pickup-and-delivery tsp.
\newblock {\em INFORMS Journal on Computing}, 32(2):263--278, 2020.

\bibitem{cplex}
{Cplex, IBM ILOG}.
\newblock {CPLEX} optimization studio 12.10, 2019.

\bibitem{bdd_ordering_fujita}
M. Fujita, Y. Lu, E. Clarke, and J. Jain.
\newblock Efficient variable ordering using abdd based sampling.
\newblock In {\em Design Automation Conference}, pages 687--692, Los Alamitos,
  CA, USA, jun 2000. IEEE Computer Society.

\bibitem{gasse2019exact_co_gnn}
Maxime Gasse, Didier Ch{\'e}telat, Nicola Ferroni, Laurent Charlin, and Andrea
  Lodi.
\newblock Exact combinatorial optimization with graph convolutional neural
  networks.
\newblock {\em arXiv preprint arXiv:1906.01629}, 2019.

\bibitem{miplib_2017}
Ambros Gleixner, Gregor Hendel, Gerald Gamrath, Tobias Achterberg, Michael
  Bastubbe, Timo Berthold, Philipp~M. Christophel, Kati Jarck, Thorsten Koch,
  Jeff Linderoth, Marco L\"ubbecke, Hans~D. Mittelmann, Derya Ozyurt, Ted~K.
  Ralphs, Domenico Salvagnin, and Yuji Shinano.
\newblock {MIPLIB 2017: Data-Driven Compilation of the 6th Mixed-Integer
  Programming Library}.
\newblock {\em Mathematical Programming Computation}, 2021.

\bibitem{globerson2008fixing}
Amir Globerson and Tommi~S Jaakkola.
\newblock Fixing max-product: Convergent message passing algorithms for {MAP}
  {LP}-relaxations.
\newblock In {\em Advances in neural information processing systems}, pages
  553--560, 2008.

\bibitem{gondzio2003parallel}
Jacek Gondzio and Robert Sarkissian.
\newblock Parallel interior-point solver for structured linear programs.
\newblock {\em Mathematical Programming}, 96(3):561--584, 2003.

\bibitem{gonzalez2020bdd}
Jaime~E Gonz{\'a}lez, Andre~A Cire, Andrea Lodi, and Louis-Martin Rousseau.
\newblock {BDD}-based optimization for the quadratic stable set problem.
\newblock {\em Discrete Optimization}, page 100610, 2020.

\bibitem{gonzalez2020integrated}
Jaime~E Gonz\'{a}lez, Andre~A Cire, Andrea Lodi, and Louis-Martin Rousseau.
\newblock Integrated integer programming and decision diagram search tree with
  an application to the maximum independent set problem.
\newblock {\em Constraints}, pages 1--24, 2020.

\bibitem{gurobi}
{Gurobi Optimization, LLC}.
\newblock {Gurobi Optimizer Reference Manual}, 2021.

\bibitem{haller2020primal}
Stefan Haller, Mangal Prakash, Lisa Hutschenreiter, Tobias Pietzsch, Carsten
  Rother, Florian Jug, Paul Swoboda, and Bogdan Savchynskyy.
\newblock A primal-dual solver for large-scale tracking-by-assignment.
\newblock In {\em AISTATS}, 2020.

\bibitem{Thrust}
Jared Hoberock and Nathan Bell.
\newblock Thrust: A parallel template library, 2010.
\newblock Version 1.7.0.

\bibitem{improved_job_sequence_bounds_hooker_2019}
John~N. Hooker.
\newblock Improved job sequencing bounds from decision diagrams.
\newblock In Thomas Schiex and Simon de Givry, editors, {\em Principles and
  Practice of Constraint Programming}, pages 268--283, Cham, 2019. Springer
  International Publishing.

\bibitem{hornakova2021making}
Andrea Hornakova, Timo Kaiser, Paul Swoboda, Michal Rolinek, Bodo Rosenhahn,
  and Roberto Henschel.
\newblock Making higher order {MOT} scalable: An efficient approximate solver
  for lifted disjoint paths.
\newblock In {\em Proceedings of the IEEE/CVF International Conference on
  Computer Vision}, pages 6330--6340, 2021.

\bibitem{huangfu2018parallelizing}
Qi Huangfu and J.~A.~J. Hall.
\newblock Parallelizing the dual revised simplex method.
\newblock {\em Math. Program. Comput.}, 10(1):119--142, 2018.

\bibitem{hutschenreiter2021fusion}
Lisa Hutschenreiter, Stefan Haller, Lorenz Feineis, Carsten Rother, Dagmar
  Kainmüller, and Bogdan Savchynskyy.
\newblock Fusion moves for graph matching.
\newblock 2021.

\bibitem{jancsary2011convergent}
Jeremy Jancsary and Gerald Matz.
\newblock Convergent decomposition solvers for tree-reweighted free energies.
\newblock In {\em Proceedings of the Fourteenth International Conference on
  Artificial Intelligence and Statistics}, pages 388--398, 2011.

\bibitem{johnson2007lagrangian}
Jason~K Johnson, Dmitry~M Malioutov, and Alan~S Willsky.
\newblock Lagrangian relaxation for {MAP} estimation in graphical models.
\newblock {\em arXiv preprint arXiv:0710.0013}, 2007.

\bibitem{kainmueller2014active}
Dagmar Kainmueller, Florian Jug, Carsten Rother, and Gene Myers.
\newblock Active graph matching for automatic joint segmentation and annotation
  of {C}. elegans.
\newblock In {\em International Conference on Medical Image Computing and
  Computer-Assisted Intervention}, pages 81--88. Springer, 2014.

\bibitem{Kappes2015}
J{\"{o}}rg~H. Kappes, Bj{\"{o}}rn Andres, Fred~A. Hamprecht, Christoph
  Schn{\"{o}}rr, Sebastian Nowozin, Dhruv Batra, Sungwoong Kim, Bernhard~X.
  Kausler, Thorben Kr{\"{o}}ger, Jan Lellmann, Nikos Komodakis, Bogdan
  Savchynskyy, and Carsten Rother.
\newblock A comparative study of modern inference techniques for structured
  discrete energy minimization problems.
\newblock {\em International Journal of Computer Vision}, 115(2):155--184,
  2015.

\bibitem{kappes_trws_pct}
Jörg~Hendrik Kappes, Markus Speth, Gerhard Reinelt, and Christoph Schnörr.
\newblock Towards efficient and exact {MAP}-inference for large scale discrete
  computer vision problems via combinatorial optimization.
\newblock In {\em 2013 IEEE Conference on Computer Vision and Pattern
  Recognition}, pages 1752--1758, 2013.

\bibitem{knuth2011art}
Donald~E Knuth.
\newblock {\em The art of computer programming, volume 4A: combinatorial
  algorithms, part 1}.
\newblock Pearson Education India, 2011.

\bibitem{kolmogorov2006convergent}
Vladimir Kolmogorov.
\newblock Convergent tree-reweighted message passing for energy minimization.
\newblock {\em IEEE transactions on pattern analysis and machine intelligence},
  28(10):1568--1583, 2006.

\bibitem{kolmogorov2014new}
Vladimir Kolmogorov.
\newblock A new look at reweighted message passing.
\newblock {\em IEEE transactions on pattern analysis and machine intelligence},
  37(5):919--930, 2014.

\bibitem{komodakis2007approximate}
Nikos Komodakis and Georgios Tziritas.
\newblock Approximate labeling via graph cuts based on linear programming.
\newblock {\em IEEE transactions on pattern analysis and machine intelligence},
  29(8):1436--1453, 2007.

\bibitem{lange2018partial}
Jan-Hendrik Lange, Andreas Karrenbauer, and Bjoern Andres.
\newblock Partial optimality and fast lower bounds for weighted correlation
  clustering.
\newblock In {\em International Conference on Machine Learning}, pages
  2892--2901. PMLR, 2018.

\bibitem{lange2021efficient}
Jan-Hendrik Lange and Paul Swoboda.
\newblock Efficient message passing for 0--1 {ILPs} with binary decision
  diagrams.
\newblock In {\em International Conference on Machine Learning}, pages
  6000--6010. PMLR, 2021.

\bibitem{lozano2018consistent}
Leonardo Lozano, David Bergman, and J~Cole Smith.
\newblock On the consistent path problem.
\newblock {\em Optimization Online e-prints}, 2018.

\bibitem{meltzer2012convergent}
Talya Meltzer, Amir Globerson, and Yair Weiss.
\newblock Convergent message passing algorithms-a unifying view.
\newblock {\em arXiv preprint arXiv:1205.2625}, 2012.

\bibitem{nair2020solving}
Vinod Nair, Sergey Bartunov, Felix Gimeno, Ingrid von Glehn, Pawel Lichocki,
  Ivan Lobov, Brendan O'Donoghue, Nicolas Sonnerat, Christian Tjandraatmadja,
  Pengming Wang, et~al.
\newblock Solving mixed integer programs using neural networks.
\newblock {\em arXiv preprint arXiv:2012.13349}, 2020.

\bibitem{cuda}
NVIDIA, Péter Vingelmann, and Frank~H.P. Fitzek.
\newblock {CUDA}, release: 11.2, 2021.

\bibitem{perumalla_gpu_mip_design_considerations}
Kalyan Perumalla and Maksudul Alam.
\newblock {\em \textnormal{Design Considerations for GPU-Based Mixed Integer
  Programming on Parallel Computing Platforms}}.
\newblock Association for Computing Machinery, New York, NY, USA, 2021.

\bibitem{ralphs2018parallel}
Ted Ralphs, Yuji Shinano, Timo Berthold, and Thorsten Koch.
\newblock Parallel solvers for mixed integer linear optimization.
\newblock In {\em Handbook of parallel constraint reasoning}, pages 283--336.
  Springer, 2018.

\bibitem{savchynskyy12efficient}
B. Savchynskyy, S. Schmidt, J{\"o}rg~H. Kappes, and Christoph Schn{\"o}rr.
\newblock Efficient {MRF} energy minimization via adaptive diminishing
  smoothing.
\newblock {\em UAI. Proceedings}, pages 746--755, 2012.
\newblock 1.

\bibitem{shekhovtsov2016solving}
Alexander Shekhovtsov, Christian Reinbacher, Gottfried Graber, and Thomas Pock.
\newblock Solving dense image matching in real-time using discrete-continuous
  optimization.
\newblock In {\em Proceedings of the 21st Computer Vision Winter Workshop
  (CVWW)}, page~13, 2016.

\bibitem{smith12_gpu_interiorpoint}
Edmund Smith, Jacek Gondzio, and Julian Hall.
\newblock {GPU} acceleration of the matrix-free interior point method.
\newblock In Roman Wyrzykowski, Jack Dongarra, Konrad Karczewski, and Jerzy
  Wa{\'{s}}niewski, editors, {\em Parallel Processing and Applied Mathematics},
  pages 681--689, Berlin, Heidelberg, 2012. Springer Berlin Heidelberg.

\bibitem{sofranac2020accelerating_domain_propagation}
Boro Sofranac, Ambros Gleixner, and Sebastian Pokutta.
\newblock Accelerating domain propagation: an efficient {GPU}-parallel
  algorithm over sparse matrices.
\newblock In {\em 2020 IEEE/ACM 10th Workshop on Irregular Applications:
  Architectures and Algorithms (IA3)}, pages 1--11. IEEE, 2020.

\bibitem{sonnerat2021learning}
Nicolas Sonnerat, Pengming Wang, Ira Ktena, Sergey Bartunov, and Vinod Nair.
\newblock Learning a large neighborhood search algorithm for mixed integer
  programs.
\newblock {\em arXiv preprint arXiv:2107.10201}, 2021.

\bibitem{swoboda2017message}
Paul Swoboda and Bjoern Andres.
\newblock A message passing algorithm for the minimum cost multicut problem.
\newblock In {\em Proceedings of the IEEE Conference on Computer Vision and
  Pattern Recognition}, pages 1617--1626, 2017.

\bibitem{swoboda2022structured_pred_archive}
Paul Swoboda, Andrea Hornakova, Paul Roetzer, and Ahmed Abbas.
\newblock Structured prediction problem archive.
\newblock {\em arXiv preprint arXiv:2202.03574}, 2022.

\bibitem{swoboda2019convex}
Paul Swoboda, Ashkan Mokarian, Christian Theobalt, Florian Bernard, et~al.
\newblock A convex relaxation for multi-graph matching.
\newblock In {\em Proceedings of the IEEE Conference on Computer Vision and
  Pattern Recognition}, pages 11156--11165, 2019.

\bibitem{swoboda2017study}
Paul Swoboda, Carsten Rother, Hassan Abu~Alhaija, Dagmar Kainmuller, and Bogdan
  Savchynskyy.
\newblock A study of lagrangean decompositions and dual ascent solvers for
  graph matching.
\newblock In {\em Proceedings of the IEEE Conference on Computer Vision and
  Pattern Recognition}, pages 1607--1616, 2017.

\bibitem{tjandraatmadja2020incorporating}
Christian Tjandraatmadja and Willem-Jan van Hoeve.
\newblock Incorporating bounds from decision diagrams into integer programming.
\newblock {\em Mathematical Programming Computation}, pages 1--32, 2020.

\bibitem{torresani2008feature}
Lorenzo Torresani, Vladimir Kolmogorov, and Carsten Rother.
\newblock Feature correspondence via graph matching: Models and global
  optimization.
\newblock In {\em European conference on computer vision}, pages 596--609.
  Springer, 2008.

\bibitem{tourani2018mplp++}
Siddharth Tourani, Alexander Shekhovtsov, Carsten Rother, and Bogdan
  Savchynskyy.
\newblock {MPLP++}: Fast, parallel dual block-coordinate ascent for dense
  graphical models.
\newblock In {\em Proceedings of the European Conference on Computer Vision
  (ECCV)}, pages 251--267, 2018.

\bibitem{tourani2020taxonomy}
Siddharth Tourani, Alexander Shekhovtsov, Carsten Rother, and Bogdan
  Savchynskyy.
\newblock Taxonomy of dual block-coordinate ascent methods for discrete energy
  minimization.
\newblock In {\em AISTATS}, 2020.

\bibitem{vineet2008cuda}
Vibhav Vineet and PJ Narayanan.
\newblock {CUDA} cuts: Fast graph cuts on the {GPU}.
\newblock In {\em 2008 IEEE Computer Society Conference on Computer Vision and
  Pattern Recognition Workshops}, pages 1--8. IEEE, 2008.

\bibitem{wang2013subproblem}
Huayan Wang and Daphne Koller.
\newblock Subproblem-tree calibration: A unified approach to max-product
  message passing.
\newblock In {\em ICML (2)}, pages 190--198, 2013.

\bibitem{werner2007linear}
Tomas Werner.
\newblock A linear programming approach to max-sum problem: A review.
\newblock {\em IEEE transactions on pattern analysis and machine intelligence},
  29(7):1165--1179, 2007.

\bibitem{werner2019relative}
Tom{\' a}{\v s} Werner, Daniel Pr{\r u}{\v s}a, and Tom{\' a}{\v s} Dlask.
\newblock Relative interior rule in block-coordinate descent.
\newblock In {\em Proceedings of the IEEE International Conference on Computer
  Vision}, 2020.
\newblock To appear.

\bibitem{wu2012chapter}
Jiadong Wu, Zhengyu He, and Bo Hong.
\newblock Chapter 5 - efficient {CUDA} algorithms for the maximum network flow
  problem.
\newblock In Wen mei W.~Hwu, editor, {\em GPU Computing Gems Jade Edition},
  Applications of GPU Computing Series, pages 55--66. Morgan Kaufmann, Boston,
  2012.

\bibitem{xu2020fast_message_passing_mrf}
Zhiwei Xu, Thalaiyasingam Ajanthan, and Richard Hartley.
\newblock Fast and differentiable message passing on pairwise markov random
  fields.
\newblock In {\em Proceedings of the Asian Conference on Computer Vision},
  2020.

\bibitem{zhang2016pairwise}
Zhen Zhang, Qinfeng Shi, Julian McAuley, Wei Wei, Yanning Zhang, and Anton Van
  Den~Hengel.
\newblock Pairwise matching through max-weight bipartite belief propagation.
\newblock In {\em Proceedings of the IEEE Conference on Computer Vision and
  Pattern Recognition}, pages 1202--1210, 2016.

\end{thebibliography}
